\documentclass{amsart}
\usepackage{amssymb}
\usepackage{tabls}

\newtheorem{proposition}{Proposition}[section]
\newtheorem{theorem}[proposition]{Theorem}
\newtheorem{lemma}[proposition]{Lemma}
\newtheorem{corollary}[proposition]{Corollary}
\theoremstyle{definition}
\newtheorem{example}[proposition]{Example}
\newtheorem{definition}[proposition]{Definition}
\newcommand{\mx}[1]{\mathrm{max}({#1})}
\newcommand{\Aut}{\mathrm{Aut}}
\newcommand{\Alg}{\mathcal{A}}
\newcommand{\dvee}{\dot\vee}
\begin{document}
\title[Representation of MV-algebras]{A representation theorem for MV-algebras}
\author{Gejza Jen\v ca}
\address 
{
Department of Mathematics\\ Faculty of Electrical Engineering and
Information Technology\\ Ilkovi\v cova~3\\ 812~19~Bratislava\\ Slovakia
}
\email{jenca@kmat.elf.stuba.sk}
\thanks{This research is supported by grant VEGA G-1/3025/06 of M\v S SR,
Slovakia and by the Science and Technology Assistance Agency under the contract
No. APVT-51-032002}
\begin{abstract}
An {\em MV-pair} is a pair $(B,G)$ where $B$ is a Boolean algebra and $G$ is a
subgroup of the automorphism group of $B$ satisfying certain conditions. 
Let $\sim_G$ be the equivalence relation on $B$ naturally associated with $G$.
We prove that for every MV-pair $(B,G)$, the effect algebra $B/\sim_G$ is an MV-effect
algebra.  Moreover, for every MV-effect algebra $M$ there is an MV-pair $(B,G)$
such that $M$ is isomorphic to $B/\sim_G$.
\end{abstract}
\maketitle
\section{Introduction}

Let $D$ be a bounded distributive lattice.
Recall, that a Boolean algebra $B(D)$ is called {\em R-generated by} $D$
iff $D$ is a $0,1$-sublattice of $B(D)$ and $D$ generates $B(D)$, as a
Boolean algebra. Given $D$, these properties determine $B(D)$ up to isomorphism.

In \cite{Jen:BARGbMVEA}, it was proved that every MV-effect algebra $M$ there
is a surjective  morphism of effect algebras $\phi_M:B(M)\to M$.  Since
$\phi_M$ is a full morphism of effect algebras, $B/\sim_{\phi_M}$ is isomorphic
to $M$. A natural question arises: is it possible to express $\phi_M$ in terms
of $B(M)$, using only the language of Boolean algebras? In this paper, we
answer this question in the affirmative. We prove that for every MV-algebra $M$
there exists a subgroup $G(M)$ of the automorphism group of $B(M)$ such that
the standard equivalence relation on $B(M)$ associated with $G(M)$ equals
$\sim_{\phi_M}$. Conversely, we give conditions under which a pair $(B,G)$
gives rise to an MV-effect algebra in aforementioned way; we call such pairs
$(B,G)$ {\em MV-pairs}.  Finally, we prove that $(B(M),G(M))$ is an MV-pair.

The origins of the main idea of this paper lie in the paper \cite{CigDubMun:AMAIfBAwaFOA}.

\section{Definitions and basic relationships}

An {\em effect algebra} is a partial algebra $(E;\oplus,0,1)$ with a binary 
partial operation $\oplus$ and two nullary operations $0,1$ satisfying
the following conditions.
\begin{enumerate}
\item[(E1)]If $a\oplus b$ is defined, then $b\oplus a$ is defined and
		$a\oplus b=b\oplus a$.
\item[(E2)]If $a\oplus b$ and $(a\oplus b)\oplus c$ are defined, then
		$b\oplus c$ and $a\oplus(b\oplus c)$ are defined and
		$(a\oplus b)\oplus c=a\oplus(b\oplus c)$.
\item[(E3)]For every $a\in E$ there is a unique $a'\in E$ such that
		$a\oplus a'=1$.
\item[(E4)]If $a\oplus 1$ exists, then $a=0$
\end{enumerate}

Effect algebras were introduced by Foulis and Bennett in their paper
\cite{FouBen:EAaUQL}. In their papers \cite{Kop:DPFS} and \cite{KopCho:DP}, K\^
opka and Chovanec introduced an essentially equivalent structure called {\em
D-poset}.  Another equivalent structure, called {\em weak orthoalgebras} was
introduced by Giuntini and Greuling in \cite{GiuGre:TaFLfUP}.  We refer to the
monograph \cite{DvuPul:NTiQS} for more information on effect algebras and
similar algebraic structures.

For brevity, we denote an effect algebra $(E;\oplus,0,1)$ by $E$.  In an effect
algebra $E$, we write $a\leq b$ iff there is $c\in E$ such that $a\oplus c=b$.
It is easy to check that every effect algebra is cancellative, thus $\leq$ is a
partial order on $E$. In this partial order, $0$ is the least and $1$ is the
greatest element of $E$.  Moreover, it is possible to introduce a new partial
operation $\ominus$; $b\ominus a$ is defined iff $a\leq b$ and then
$a\oplus(b\ominus a)=b$.  It can be proved that $a\oplus b$ is defined iff
$a\leq b'$ iff $b\leq a'$. Therefore, we denote the domain of $\oplus$
by $\perp$.

Let $E_1,E_2$ be effect algebras. A mapping $\phi:E_1\mapsto E_2$ is called a
{\em morphism of effect algebras} iff $\phi(1)=1$ and for all $a,b\in E$, the
existence of $a\oplus b$ implies the existence of $\phi(a)\oplus\phi(b)$ and
$\phi(a\oplus b)=\phi(a)\oplus\phi(b)$.  A morphism $\phi:E_1\to E_2$ is {\em
full} iff whenever $\phi(a)\perp\phi(b)$ and
$\phi(a)\oplus\phi(b)\in\phi(E_1)$, then there are $a_1,b_1\in E_1$ such that
$\phi(a)=\phi(a_1)$, $\phi(b)=\phi(b_1)$ and $a_1\perp b_1$.  A morphism $\phi$
is an {\em isomorphism} iff $\phi$ is bijective and full. Note that even if
both $E_1$ and $E_2$ are lattice ordered, a morphism of effect algebras need
not to preserve joins and meets. 

An {\em MV-algebra} (c.f. \cite{Cha:AAoMVL}, \cite{Mun:IoAFCSAiLSC}) is a
$(2,1,0)$-type algebra $(M;\boxplus,\lnot,0)$, 
such that $\boxplus$
satisfying the
identities $(x\boxplus y)\boxplus z=x\boxplus (y\boxplus z)$,
$x\boxplus z=y\boxplus x$, $x\boxplus0=x$, $\lnot\lnot x=x$, $x\boxplus\lnot 0=\lnot 0$ and
$$
x\boxplus\lnot(x\boxplus\lnot y)=y\boxplus\lnot(y\boxplus\lnot x)\text{.}
$$
On every MV-algebra, a partial order $\leq$ is defined
by the rule
$$
x\leq y\Longleftrightarrow y=x\boxplus\lnot(x\boxplus\lnot y).
$$
In this partial order, every MV-algebra is a distributive lattice bounded by
$0$ and $\lnot 0$.

An {\em MV-effect algebra} is a lattice ordered effect algebra $M$ in which,
for all $a,b\in M$, $(a\lor b)\ominus a=b\ominus (a\land b)$. It is proved
in \cite{ChoKop:BDP} that there is a natural, one-to one correspondence between
MV-effect algebras and MV-algebras given by the following rules.
Let $(M,\oplus,0,1)$ be an MV-effect algebra. Let $\boxplus$ 
be a total operation given by $x \boxplus y=x\oplus(x'\land y)$. Then
$(M,\boxplus,',0)$ is an MV-algebra. Similarly, let 
$(M,\boxplus,\lnot,0)$ be an MV-algebra. Restrict the operation
$\boxplus$ to the pairs $(x,y)$ satisfying $x\leq y'$ and call the
new partial operation $\oplus$. Then $(M,\oplus,0,\lnot 0)$ is an MV-effect algebra.

Among lattice ordered effect algebras, MV-effect algebras can be characterized
in a variety of ways. Three of them are given in the following
proposition.
\begin{proposition}
\cite{BenFou:PSEA}, \cite{ChoKop:BDP}
Let $E$ be a lattice ordered effect algebra. The following are equivalent
\begin{enumerate}
\item[(a)] $E$ is an MV-effect algebra.
\item[(b)] For all $a,b\in E$, $a\land b=0$ implies $a\leq b'$.
\item[(c)] For all $a,b\in E$, $a\ominus(a\land b)\leq b'$.
\item[(d)] For all $a,b\in E$, there exist $a_1,b_1,c\in E$ such that
$a_1\oplus b_1\oplus c$ exists, $a_1\oplus c=a$ and $b_1\oplus c=b$.
\end{enumerate}
\end{proposition}

\noindent{\bf Notation.} In what follows, we will deal with an MV-effect
algebra $M$ and a Boolean algebra $B(M)$ such that $M$ is a 0,1-sublattice of
$B(M)$.  In this particular situation, a small notational problem arises: both
$M$ and $B(M)$ are MV-effect algebras, but the $\oplus,\ominus$ and $~'$
operations on $B(M)$ and $M$ differ.  To avoid confusion, we denote the partial
operation of disjoint join (the $\oplus$ of Boolean algebras) on a Boolean
algebra by $\dvee$. The partial difference of comparable elements and the
complement in a Boolean algebra are denoted by $\setminus$ and $~^\complement$,
respectively.

Let $D$ be a bounded distributive lattice. Up to isomorphism, there exists a unique Boolean algebra $B(D)$ such that $D$ is a $0,1$-sublattice of $B(D)$ and
$B$ generates $B(D)$ as a Boolean algebra. 
This Boolean algebra is called the Boolean algebra R-generated by $D$. 
We refer to \cite{Gra:GLT}, section II.4, for an overview of results
concerning R-generated Boolean algebras. See also \cite{Has:ITfL} and 
\cite{MacNei:EoaDLtaBR}.
For every element $x$ of $B(D)$, there exists a finite chain 
$x_1\leq\ldots\leq x_n$ in $D$ such that $x=x_1+\ldots+x_n$. Here, $+$ denotes
the symmetric difference, as in Boolean rings. We then say than $\{x_i\}_{i=1}^n$ is
a {\em $D$-chain representation} of $x$. It is easy to see that every
element of $B(D)$ has a $D$-chain representation of even length.
Note that, for $n=2k$ we have
$$
x=x_1+\dots+x_{2k}=(x_{2k}\setminus x_{2k-1})
	\oplus\dots\oplus(x_2\setminus x_1).
$$

If $D_1,D_2$ are bounded distributive lattices and $\psi:D_1\to D_2$ is
a $0,1$-lattice homomorphism, then $\psi$ uniquely extends to
a homomorphism of Boolean algebras $\psi^*:B(D_1)\to B(D_2)$.
Similarly, if $[0,a]_D$ is an interval in a bounded distributive
lattice $D$, then $B([0,a]_D)$ is naturally isomorphic to
the interval $[0,a]_{B(D)}$.

\begin{theorem}\cite{Jen:BARGbMVEA}
\label{thm:main}
Let $M$ be an MV-effect algebra. The mapping $\phi_M:B(M)\to M$
given by
$$
\phi_M(x)=\bigoplus_{i=1}^n(x_{2i}\ominus x_{2i-1}),
$$
where $\{x_i\}_{i=1}^{2n}$ is a $M$-chain representation of $x$,
is a surjective morphism of effect algebras.
\end{theorem}

We note that the value of $\phi_M(x)$ does not depend on the choice
of the $M$-chain representation of $x$. Obviously, for all
$x\in M$, $\{x,0\}$ is a $M$-chain representation of $x$. Therefore,
$\phi_M(x)=x\ominus 0=x$, so every $x\in M$ is a fixpoint
of $\phi_M$.

\begin{example}
\label{ex:chain}
Let $M$ be an MV-effect algebra, which is totally ordered.
By \cite{Gra:GLT}, Corollary II.4.19, $B(M)$ is
isomorphic to the Boolean algebra of all subsets of
$M$ of the form $[a_1,b_1)\dot\cup\ldots\dot\cup[a_n,b_n)$.
Here, we denote $[a,b)=\{x\in M:a\leq x<b\}$.
The $\phi_M:B(M)\to M$ morphism is then given by
$$
\phi_M([a_1,b_1)\dot\cup\ldots\dot\cup[a_n,b_n))=
(b_1\ominus a_1)\oplus\ldots\oplus(b_n\ominus a_n).
$$
\end{example}

\begin{example}
In this example, $[0,1]$ denotes the closed real unit interval.
Let $C_{[0,1]}$ be the MV-effect algebra of all real continuous functions
$f:[0,1]\to[0,1]$. Let $B$ be the Boolean algebra
$$
\prod_{x\in [0,1]}B([0,1]),
$$
where $B([0,1])$ is the Boolean algebra generated by semiopen intervals 
as described in Example \ref{ex:chain}.
It is obvious that $C_{[0,1]}$, as a bounded lattice, can be embedded
into $B$ by a mapping $\gamma:E\to B$ given by
$\gamma(f)=\bigl(\bigl[0,f(x)\bigr)\bigr)_{x\in[0,1]}$. The image
of $E$ under $\gamma$ then generates a Boolean subalgebra of
$B$, which we can identify with $B(C_{[0,1]})$.
The $\phi_{C_{[0,1]}}:B(C_{[0,1]})\to C_{[0,1]}$ mapping can then be 
constructed as follows.

Let $(A_x)_{x\in[0,1]}\in B(C_{[0,1]})$.
Fix $x\in[0,1]$ and write $A_x=[a_1,b_1)\dot\cup\ldots\dot\cup[a_n,b_n)$.
The value of the continuous function $\phi_{C_{[0,1]}}((A_x)_{x\in[0,1]})$
at $x$ is then equal to $(b_1\ominus a_1)\oplus\ldots\oplus(b_n\ominus a_n)$.
\end{example}

Let $E$ be an effect algebra. A relation $\sim$ on $E$ is 
a {\em weak congruence} iff the following conditions are satisfied.
\maketitle
\begin{enumerate}
\item[(C1)]$\sim$ is an equivalence relation.
\item[(C2)]If $a_1\sim a_2$, $b_1\sim b_2$ and  $a_1\oplus b_1,a_2\oplus b_2$ exist, then
$a_1\oplus b_1\sim a_2\oplus b_2$.
\end{enumerate}

If $E$ is an effect algebra and $\sim$ is a weak congruence on $E$, the
quotient $E/\sim$ ($\oplus$ is defined on $E/\sim$ in an obvious way) 
need not to be a partial abelian monoid, since the associativity condition may 
fail (c.f. \cite{GudPul:QoPAM}). This fact motivates the study of 
sufficient conditions for a weak congruence
to preserve associativity. The following condition was considered in 
\cite{ChePul:SILiPAM}.
\begin{enumerate}
\item[(C5)]If $a\sim b\oplus c$, then there are $b_1,c_1$ such that $b_1\sim b$, $c_1\sim c$, $b_1\oplus c_1$ exists and $a=b_1\oplus c_1$.
\end{enumerate}

In \cite{ChePul:SILiPAM}, it was proved that for a partial abelian monoid $P$ and
a weak congruence $\sim$, satisfying (C5), the quotient $P/\sim$ is again a
partial abelian monoid. Moreover, it is easy to prove that the eventual 
positivity of $P$ is preserved for such $\sim$.
However, for an effect algebra $E$, the (C5) property of $\sim$ does not
guarantee that the ${}'$ operation is preserved by $\sim$.
If ${}'$ is preserved by $\sim$, that means, if condition
\begin{enumerate}
\item[(C6)]If $a\sim b$, then $a'\sim b'$.
\end{enumerate}
is satisfied, then $E/\sim$ is an effect algebra. A relation on an effect algebra satisfying (C1),(C2),(C5),(C6) is called {\em an effect algebra congruence}.
For every effect algebra congruence $\sim$ on an effect algebra $E$, the
mapping $a\to [a]_\sim$ is a full morphism of effect algebras.

We refer the interested reader to \cite{Pul:CiPAS} and \cite{GudPul:QoPAM} for further 
details concerning congruences on effect algebras and partial abelian monoids.

The (b) and (c) of the following lemma are just two equivalent $\perp$-to-$\leq$
reformulations of the (C3) property from \cite{GudPul:QoPAM}. Thus, the lemma
is (implicitly) well known, but we cannot find it in print. 
\begin{lemma}
\label{lemma:select}
Let $\sim$ be a congruence on an effect algebra $E$. For all $x,y\in E$, the
following are equivalent.
\begin{enumerate}
\item[(a)]$[x]_\sim\leq[y]_\sim$.
\item[(b)]There is $x_1\sim x$ such that $x_1\leq y$.
\item[(c)]There is $y_1\sim y$ such that $x\leq y_1$.
\end{enumerate}
\end{lemma}
\begin{proof}~

(b)$\implies$(a) and (c)$\implies$(a) are trivial.

(a)$\implies$(b): As $[x]_\sim\leq[y]_\sim$, there is $u\in E$ such that 
$[x]_\sim\oplus[u]_\sim=[y]_\sim$. This implies that there are
$x_0,u_0\in E$ such that $x_0\sim x$, $u_0\sim u$,
$x_0\oplus u_0$ exists, and $x_0\oplus u_0\sim y$.
By the (C5) property, there are $x_1,u_1$ such that
$x_1\sim x_0$, $u_1\sim u_0$, $x_1\oplus u_1$ exists,
and $x_1\oplus u_1=y$.

(a)$\implies$(c): By the (C6) property, $[y']_\sim\leq [x']_\sim$.
As (a)$\implies$(b), there is $z\sim y'$ such that $z\leq x'$
and this is equivalent with $x\leq z'$. By the (C6) property,
$z\sim y'$ iff $z'\sim y$ and we can put $y_1=z'$.
\end{proof}

Recall that an effect algebra $E$ satisfies the {\em Riesz decomposition
property} iff for all $u,v_1,v_2\in E$, $u\leq v_1\oplus v_2$ iff there are
$u_1,u_2$ such that $u_1\leq v_1$, $u_2\leq v_2$ and $u=u_1\oplus u_2$.  A
lattice ordered effect algebra is an MV-effect algebra iff it satisfies the
Riesz decomposition property. There are non-lattice ordered effect algebras
satisfying the Riesz decomposition property, for example the effect algebra of
all polynomial functions $[0,1]_{\mathbb R}\to[0,1]_{\mathbb R}$. By
\cite{Rav:OaSToEA}, every
effect algebra satisfying the Riesz decomposition property can be
embedded, as an interval in the positive cone, into a partially ordered abelian
group satisfying the Riesz decomposition property. This result is a generalization of the famous result by Mundici from \cite{Mun:IoAFCSAiLSC}.

An effect algebra satisfies the {\em Riesz interpolation property} iff for all
elements $u_1,u_2,v_1,v_2$ such that $u_i\leq v_j$ for all $i,j\in\{1,2\}$,
there is an element $x$ such that $x$ is an upper bound of $u_1,u_2$ and a
lower bound of $v_1,v_2$.  If an effect algebra satisfies the Riesz
decomposition property, then it satisfies the Riesz interpolation property. The
opposite implication is not true, since every lattice ordered effect algebra
satisfies the Riesz interpolation property, but there exist (obviously) some
effect algebras that are lattice ordered and non-MV.

\section{From MV-pairs to MV-effect algebras}

Let $B$ be a Boolean algebra. We write $\Aut(B)$ for the group 
of all automorphisms of $B$. Let $G$ be a subgroup of $\Aut(B)$.
For $a,b\in B$, we write $a\sim_G b$ iff there exists $f\in G$ such that
$b=f(a)$. Obviously, $\sim_G$ is an equivalence relation. We write
$[a]_G$ for the equivalence class of an element $a$ of $B$.

A pair $(B,G)$, where $B$ is a Boolean algebra and $G$ is a subgroup of
$\Aut(B)$ is called {\em a BG-pair}. BG-pairs are a well-established topic
in the theory of Boolean algebras, see for example Chapter 15 of the
handbook \cite{Kop:HoBA}.

Let $(P,\leq)$ be a poset. Let us write,
$$
\mx{P}=\{m\in P: x\leq m\implies x=m\},
$$
that means, $\mx{P}$ is the set of all maximal elements of the poset $P$.

Let $B$ be a Boolean algebra, let $G$ be a subgroup of $\Aut(B)$.
For all $a,b\in B$, we write
\begin{align*}
L(a,b)&=\{a\wedge f(b):f\in G\}\text{ and}\\
L^+(a,b)&=\{g(a)\wedge f(b): f,g\in G\}.
\end{align*}
Note that $L(a,b)\subseteq L^+(a,b)$ and that $L^+(a,b)$ is closed with respect
to any $h\in G$; this implies that $L^+(a,b)$ is a union of equivalence
classes of $\sim_G$.

\begin{definition}
Let $B$ be a Boolean algebra, let $G$ be a subgroup of $\Aut(B)$.
We say that $(B,G)$ is an {\em MV-pair} iff the following two conditions are
satisfied.
\begin{enumerate}
\item[(MVP1)] For all $a,b\in B$, $f\in G$ such that $a\leq b$ and $f(a)\leq b$,
there is $h\in G$ such that $h(a)=f(a)$ and $h(b)=b$.
\item[(MVP2)] For all $a,b\in B$ and $x\in L(a,b)$, there exists
$m\in\mx{L(a,b)}$ with $m\geq x$.
\end{enumerate}
\end{definition}
\begin{example}
For every finite Boolean algebra $B$, $(B,\Aut(B))$ is an MV-pair.
\end{example}

\begin{example}
Let $B$ be a Boolean algebra with three atoms $a_1,a_2,a_3$.
The mapping $f$ given by
$$
\begin{array}{|r||c|c|c|c|c|c|c|c|}
\hline
x&0&a_1&a_2&a_3&a_1^\complement&a_2^\complement&a_3^\complement&1\\
\hline
f(x)&0&a_2&a_3&a_1&a_2^\complement&a_3^\complement&a_1^\complement&1\\
\hline
\end{array}
$$
is an automorphism of $B$ and $G=\{id,f,f^2\}$ is a subgroup
of $\Aut(B)$. However, $(B,G)$ is not an MV-pair. Indeed,
we have $a_1\leq a_3^\complement$ and $f(a_1)=a_2\leq a_3^\complement$,
but there is no $h\in G$ such that $h(a_1)=f(a_1)$ and $h(a_3^\complement)=a_3^\complement$.
\end{example}
\begin{example}
Let $B$ be the Boolean algebra of all Borel subsets of the real unit interval
$[0,1]_{\mathbb R}$ that are unions of a finite number of intervals.
(as usual, we identify the Borel sets that differ by a set
of measure $0$.) Let $W$ the subgroup of the permutation group of
$[0,1]_{\mathbb R}$ that is generated by the set of all bijections $p_{a,b}$ given by
$$
p_{a,b}(x)=
\begin{cases}
x&\text{if $x\in [0,a]$,}\\
a+b-x\;&\text{if $x\in (a,b)$,}\\
x&\text{if $x\in[b,1]$,}
\end{cases}
$$
where $0\leq a\leq b\leq 1$. For every $p\in W$, let $f_p$ be
the mapping $f_p:B\to B$ given by
$f_p(X)=p(X)$ and let $G=\{f_p:p\in W\}$. Obviously,
$G$ is a subgroup of $\Aut(B)$. Then $(B,G)$ is an MV-pair; the
proof of this fact is a bit longer, but straightforward.
Note that every $f_p\in G$ preserves measure.
\end{example}

\begin{example}
Let $2^{\mathbb Z}$ be the Boolean algebra of all subsets of $\mathbb Z$.
Then $(2^{\mathbb Z},\Aut(2^{\mathbb Z}))$ is not an MV-pair. Indeed, let $f\in\Aut(2^{\mathbb Z})$ be the
automorphism of $2^{\mathbb Z}$ associated with the permutation
$f(n)=n+1$.
Let $A=B=\mathbb N$. We see that
$f(A)=A\setminus\{0\}$, $A\subseteq B$ and $f(A)\subseteq B$.
However,  
there is no $h\in\Aut(2^{\mathbb Z})$ such that $h(A)=f(A)$ and $h(B)=B$, simply because
$A=B$ implies that $h(A)=h(B)$, but $f(A)\neq B$.
\end{example}

The (MVP1) condition can be reformulated:
\begin{lemma}
\label{lemma:charmvp1}
Let $B$ be a Boolean algebra, let $G$ be a subgroup of $\Aut(B)$.
Then the following conditions are equivalent.
\begin{enumerate}
\item[(a)] (MVP1)
\item[(b)] For all $a,b\in B$, $f\in G$ such that $a\leq b$ and $a\leq f(b)$,
there is $h\in G$ such that $h(b)=f(b)$ and $h(a)=a$.
\item[(c)] For all $a,b\in B$, $f\in G$ such that $a\wedge b=0$ and $a\wedge
f(b)=0$, there is $h\in G$ such that $h(b)=f(b)$ and $h(a)=a$.
\end{enumerate}
\end{lemma}
\begin{proof}~

(a)$\implies$(b): Replace $a$ with $b^\complement$ and $b$ with $a^\complement$
and apply the fact that $f$ is an automorphism.

(b)$\implies$(c): Replace $b$ with $b^\complement$.

(c)$\implies$(a): Replace $b$ with $a$ and $a$ with $b^\complement$.
\end{proof}

\begin{lemma}
\label{lemma:maxtomax}
Let $(B,G)$ be an MV-pair, let $a,b\in B$ and let $m$ be a maximal element of
$L(a,b)$. For all $f\in G$, $f(m)$ is a maximal element of
$L^+(a,b)$.
\end{lemma}
\begin{proof}
Suppose that there is some element in $y\in L^+(a,b)$ with
$y\geq f(m)$ and write $y=g_1(a)\wedge f_1(b)$, where $g_1,f_1\in G$.
Since $m\in L(a,b)$, $a\geq m$ and since 
$$
a\wedge g_1^{-1}\bigl(f_1(b)\bigr)=g_1^{-1}\bigl(g_1(a)\wedge f_1(b)\bigr)=g_1^{-1}(y)\geq
g_1^{-1}\bigl(f(m)\bigr)=(g_1^{-1}\circ f)(m),
$$
we see that $a\geq (g_1^{-1}\circ f)(m)$.

By (MVP1), $a\geq(g_1^{-1}\circ f)(m)$ and $a\geq m$ imply that there exists $h\in G$ such that $h(a)=a$ and 
$h(m)=(g_1^{-1}\circ f)(m)$.
We apply $h^{-1}$ to both sides of the inequality
$$
a\wedge g_1^{-1}\bigl(f_1(b)\bigr)\geq(g_1^{-1}\circ f)(m),
$$
to obtain
$$
h^{-1}\Bigl(a\wedge g_1^{-1}\bigl(f_1(b)\bigr)\Bigr)=
a\wedge h^{-1}\Bigl(g_1^{-1}\bigl(f_1(b)\bigr)\Bigr)\geq
h^{-1}\bigl((g_1^{-1}\circ f)(m)\bigl)=m
$$
Since $m$ is a maximal element of $L(a,b)$, $a\wedge
h^{-1}\Bigl(g_1^{-1}\bigl(f_1(b)\bigr)\Bigr)\geq m$ implies that $a\wedge
h^{-1}\Bigl(g_1^{-1}\bigl(f_1(b)\bigr)\Bigr)=m$. After we apply the mapping
$g_1\circ h$ on both sides of the latter equality we obtain $y=g_1(a)\wedge
f_1(b)=f(m)$. Thus, $f(m)$ is maximal in $L^+(a,b)$.
\end{proof}
Note that Lemma \ref{lemma:maxtomax} implies that
$\mx{L(a,b)}\subseteq\mx{L^+(a,b)}$.
\begin{corollary}
Let $(B,G)$ be an MV-pair.  For all $a,b\in B$ and $x\in L^+(a,b)$, there exists
$m\in\mx{L^+(a,b)}$ with $m\geq x$.
\end{corollary}
\begin{proof}
As $x\in L^+(a,b)$, we have $x=g_1(a)\wedge f_1(b)$ for some $f_1,g_1\in G$. Then
$$
g_1^{-1}\bigl(g_1(a)\wedge f_1(b)\bigr)=a\wedge g_1^{-1}\bigl(a_1(b)\bigr)\in L(a,b).
$$
By (MVP2), there is $m\in\mx{L(a,b)}$ with $m\geq a\wedge
g_1^{-1}\bigl(a_1(b)\bigr)$.
This implies that $g_1(m)\geq g_1(a)\wedge f_1(b)$.
By Lemma \ref{lemma:maxtomax}, $g_1(m)\in\mx{L^+(a,b)}$.
\end{proof}
\begin{theorem}
\label{thm:MVPtoMV}
Let $(B,G)$ be an MV-pair. Then $\sim_G$ is an effect algebra congruence on $B$ and
$B/\sim_G$ is an MV-effect algebra.
\end{theorem}
\begin{proof}
We shall prove that the equivalence $\sim_G$ is an effect congruence.
It is easy to see that $\sim_G$ preserves the $~^\complement$ operation, so (C6) is
satisfied. To prove (C5), let $a_1,a_2\in B$ be such that $a_1\dvee a_2$ exists and
$a_1\dvee a_2\sim_G b$. Then there is $f\in G$ such that
$f(a_1\dvee a_2)=b$ and we may put $b_1=f(a_1)$ and $b_2=f(a_2)$.

Let us prove (C2).
Let $a_1,a_2,b_1,b_2\in B$ be such that $a_1\sim_G a_2$, $b_1\sim_G b_2$,
and $a_1\dvee b_1,a_2\dvee b_2$ exist. There are $f_a,f_b\in G$ such that
$f_a(a_1)=a_2$ and $f_b(b_1)=b_2$. 

We see that $b_2^\complement\geq a_2$ and that implies 
$$
b_1^\complement=f_b^{-1}(b_2^\complement)\geq f_b^{-1}(a_2)=f_b^{-1}\bigl(f_a(a_1)\bigr)=(f_b^{-1}\circ
f_a)(a_1).
$$
By (MVP1), $a_1\leq b_1^\complement$ and $(f_b^{-1}\circ f_a)(a_1)\leq b_1^\complement$ imply that
there is $h\in G$ such that $h(a_1)=(f_b^{-1}\circ f_a)(a_1)$ and
$h(b_1^\complement)=b_1^\complement$. Therefore,
$$
f_b\bigl(h(a_1\dvee b_1)\bigr)=f_b\bigl(h(a_1)\dvee h(b_1)\bigr)=
f_b\bigl((f_b^{-1}\circ f_a)(a_1)\dvee b_1\bigr)=
f_a(a_1)\dvee f_b(b_1)=a_2\dvee b_2,$$ 
and $a_1\dvee b_1\sim_G a_2\dvee b_2$.

Since $\sim_G$ is an effect congruence, $B/\sim_G$ is an effect algebra.
By Proposition 4.3 of \cite{JenPul:QoPAMatRDP}, since $B$ satisfies the Riesz decomposition property, $B/\sim_G$
satisfies the Riesz decomposition property as well. It remains to prove
that $B/\sim_G$ is a lattice. Since an effect algebra is a lattice iff
it is a (join or meet) 
semilattice, it suffices to prove that for all $a,b\in B$,
$[a]_G\wedge [b]_G$ exists in $B/\sim_G$.

Let $a,b\in B$. We shall prove that every common lower bound of $[a]_G,[b]_G$
is under a maximal common lower bound of $[a]_G,[b]_G$.

If $[c]_G\leq[a]_G,[b]_G$ then, by Lemma \ref{lemma:select},
there is
$c_1\sim_G c$ such that $c_1\leq a$ and, again by Lemma \ref{lemma:select},
$b_1\sim_G b$ such that $c_1\leq b$. As $b_1\sim_G b$, there is
$f\in G$ such that $b_1=f(b)$. Thus,
$$
c\sim_G c_1\leq a\wedge f(b)\in L(a,b).
$$
By (MVP2), there is $m\in\mx{L(a,b)}$ with $a\wedge f(b)\leq m$.
Obviously, $m\in L(a,b)$ implies that $[m]_G\leq[a]_G,[b]_G$.
Therefore, for every common lower bound $[c]_G$
of $[a]_G,[b]_G$, there is $m\in\mx{L(a,b)}$ such that
$$
[c]_G\leq[m]_G\leq[a]_G,[b]_G.
$$
Let us prove that $[m]_G$ is a maximal common lower bound
of $[a]_G,[b]_G$ in $B/\sim_G$. Suppose that 
$$
[m]_G\leq[x]_G\leq[a]_G,[b]_G.
$$
By Lemma \ref{lemma:select}, there are $m_1\sim_G m$, 
$x_1\sim_G x$ and $b_1\sim_G b$ such that
$$
m_1\leq x_1\leq a,b_1.
$$
There is $f\in G$ such that $b_1=f(b)$.
We see that $x_1\leq a\wedge f(b)\in L(a,b)\subseteq L^+(a,b)$.
There is $g\in G$ such that $m_1=g(m)$. By Lemma \ref{lemma:maxtomax},
$m_1=g(m)$ is maximal element of $L^+(a,b)$. Therefore, $m_1=a\wedge f(b)$ and
hence $x_1=m_1$. This implies that $[m]_G=[x]_G$.

Let $[m_1]_G,[m_2]_G$ be maximal common lower bounds of $[a]_G,[b]_G$.
Since $B/\sim_G$ satisfies the Riesz decomposition property,
$B/\sim_G$ satisfies the Riesz interpolation property.
By the Riesz interpolation property, there is $[m]_G$ such that
$[m_1]_G,[m_2]_G\leq[m]_G\leq[a]_G,[b]_G$.
Since $[m_1]_G,[m_2]_G$ are maximal, $[m_1]_G=[m]_G=[m_2]_G$.
Since every common lower bound of $[a]_G,[b]_G$ is under a maximal one,
and there is a single maximal common lower bound of $[a]_G,[b]_G$,
$[a]_G\wedge [b]_G$ exists.

Note that we have proved that
$[a]_G\wedge[b]_G=L^+(a,b)$. In particular, $L^+(a,b)$ is a single
equivalence class of $\sim_G$.
\end{proof}
In what follows we shall denote the MV-effect algebra
arising from an MV-pair $(B,G)$ in the way indicated above 
by $\Alg(B,G)$.

\section{From MV-effect algebras to MV-pairs}
We have proved that for every MV-pair $(B,G)$ there is an MV-effect algebra
$\Alg(B,G)$ arising from it. In this section, we shall prove
that for every MV-effect algebra there is a MV-pair $(B,G)$ such that
$\Alg(B,G)\simeq M$.

Let $M$ be an MV-effect algebra. Let $S$ be a subset of $B(M)$. We say that a mapping
$f:S\to B$ is {\em $\phi_M$-preserving} iff, for all $x\in S$,
$\phi_M(x)=\phi_M(f(x))$ or, in other words, $\phi_M$ restricted to $S$ equals
$\phi_M\circ f$.

\begin{theorem}
\label{thm:MVtoMVP}
Let $M$ be an MV-effect algebra. Let $G(M)$ be the set of all
$\phi_M$-preserving automorphisms of $B(M)$.
Then $(B(M),G(M))$ is an MV-pair and
$\Alg(B(M),G(M))$ is isomorphic to $M$.
\end{theorem}
We have divided the proof into a sequence of lemmas.  In this section, $M$ is
an MV-effect algebra and $G(M)$ is the subgroup of $\Aut(B(M))$ described in
Theorem \ref{thm:MVtoMVP}.
\begin{lemma}
\label{lemma:starter}
Let $c,d\in M$, $d\leq c$.
There is a $\phi_M$-preserving isomorphism
$$
\psi:B([0,c\ominus d]_M)\to [0,c\setminus d]_{B(M)}
$$
\end{lemma}
\begin{proof}
Consider the mapping
$
\psi_0:[0,c\ominus d]_M\to[0,c\setminus d]_{B(M)},
$
given by $\psi_0(x)=(x\oplus d)\setminus d$.
We see that $\psi_0(0)=0$, $\psi_0(c\ominus d)=c\setminus d$ and,
since $\psi_0$ is just a composition of a translation in $M$ and a translation
in $B(M)$, $\psi_0$ preserves joins and meets. Moreover, it is easy to see
that $\psi_0$ is injective, hence $\psi_0$ is a $0,1$-lattice embedding
of $[0,c\ominus d]_M$ into $[0,c\setminus d]_{B(M)}$. We shall prove that
the range of $\psi_0$ R-generates the Boolean algebra $[0,c\setminus
d]_{B(M)}$. $\psi_0$ then uniquely extends to an isomorphism
$\psi:B([0,c\ominus d]_M)\to [0,c\setminus d]_{B(M)}$.

Let $x\in[0,c\setminus d]_{B(M)}$. Let $\{x_i\}_{i=1}^{2n}$ be an
$M$-chain representation of $x$. For all $1\leq i\leq n$,
$x_{2i}\setminus x_{2i-1}\leq c\setminus d$. By elementary Boolean calculus,
this implies that
$$
x_{2i}\setminus x_{2i-1}=\bigl((x_{2i}\vee d)\wedge c\bigr)\setminus
			\bigl((x_{2i-1}\vee d)\wedge c\bigr).
$$

For all $1\leq j\leq 2n$, $(x_j\vee d)\wedge c\in[d,c]$
Therefore, $x$ has a $M$-chain representation
$\{y_j\}_{j=1}^{2n}\subseteq[d,c]_M$. Since, for all $1\leq i\leq n$,
$$
y_{2i}\setminus y_{2i-1}=(y_{2i}\setminus d)\setminus(y_{2i-1}\setminus d),
$$
$\{y_i\setminus d\}_{i=1}^{2n}$ is a chain representation of $x$.
It remains to observe that, for all $1\leq i\leq 2n$,
$$
y_i\setminus d=\bigl((y_i\ominus d)\oplus d\bigr)\setminus d=\phi_0(y_i\ominus d)
$$
and that $y_i\ominus d\in[0,c\ominus d]_M$.
Thus, every element of $[0,c\setminus d]_{B(M)}$ has a $\psi_0([0,c\ominus
d]_M)$-chain representation.

Let us prove that $\psi$ is a $\phi_M$-preserving mapping.
Let $z\in B([0,c\ominus d]_M)$, let $\{z_i\}_{i=1}^{2n}$ be a $[0,c\ominus
d]_M$-chain representation of $z$. Then
\begin{align*}
\phi_M\bigl(\psi(z)\bigr)=\phi_M\bigl(\psi(\dot\vee_{i=1}^n(z_{2i}\setminus
	z_{2i-1}))\bigr)=\\
	=\phi_M\bigl(\dot\vee_{i=1}^n\psi(z_{2i}\setminus
	z_{2i-1})\bigr)
	=\bigoplus_{i=1}^n\phi_M\bigl(\psi(z_{2i}\setminus
	z_{2i-1})\bigr)
\end{align*}
and, for all $1\leq i\leq n$,
\begin{align*}
\phi_M\bigl(\psi(z_{2i}\setminus z_{2i-1})\bigr)=
\phi_M\bigl(\psi(z_{2i})\setminus\psi(z_{2i-1})\bigr)=\\
=\phi_M\bigl(((z_{2i}\oplus d)\setminus d)\setminus((z_{2i}\oplus d)\setminus
d)\bigr)=\\
=\phi_M\bigl((z_{2i}\oplus d)\setminus(z_{2i}\oplus d)\bigr)=
\phi_M(z_{2i}\oplus d)\ominus\phi_M(z_{2i}\oplus d)=\\
=(z_{2i}\oplus d)\ominus (z_{2i-1}\oplus d)=z_{2i}\ominus z_{2i-1}=\phi_M(z_{2i}\setminus z_{2i-1}).
\end{align*}
so we obtain
$$
\phi_M(\psi(z))=\bigoplus_{i=1}^n\phi_M\bigl(\psi(z_{2i}\setminus
	z_{2i-1})\bigr)=\bigoplus_{i=1}^n\phi_M(z_{2i}\setminus z_{2i-1})=
	\phi_M(z).
$$
\end{proof}
\begin{corollary}
\label{coro:twice}
Let $c_1,d_1,c_2,d_2\in M$ be such that
$c_1\geq d_1$, $c_2\geq d_2$ and $c_1\ominus d_1=c_2\ominus d_2$.
There is a $\phi_M$-preserving isomorphism
$\psi:[0,c_1\setminus d_1]_{B(M)}\to[0,c_2\setminus d_2]_{B(M)}$.
\end{corollary}
\begin{proof}
Use Lemma \ref{lemma:starter} twice.
\end{proof}
\begin{lemma}
\label{lemma:continuer}
For every $a\in B(M)$, there is a $\phi_M$-preserving isomorphism
of Boolean algebras $\psi:B([0,\phi_M(a)]_M)\to[0,a]_{B(M)}$.
\end{lemma}
\begin{proof}
Let $\{a_i\}_{i=1}^{2n}$ be an $M$-chain representation of $a$. Then
$\{a_{2i}\setminus a_{2i-1}\}_{i=1}^n$ is a decomposition of unit in the
Boolean algebra $[0,a]_{B(M)}$ and $\phi_M(a)=\bigoplus_{i=1}^n(a_{2i}\ominus
a_{2i-1}$.
For $j\in\{0,\dots,n\}$, write $b_j=\bigoplus_{i=1}^j(a_{2i}\ominus a_{2i-1})$.
Then $\{b_j\}_{j=0}^n$ is a finite chain in $[0,\phi_M(a)]_M$ with
$b_0=0$ and $b_n=\phi_M(a)$. Thus,
$\{b_j\setminus b_{j-1}\}_{j=1}^n$ is a decomposition of unit in the Boolean
algebra
$B([0,\phi_M(a)]_M)$. For every $x\in B([0,\phi_M(a)]_M)$,
$x={\dot\bigvee}_{j=1}^n x\wedge(b_j\setminus b_{j-1})$.
Since, for all $j$, $b_j\ominus b_{j-1}=a_{2j}\ominus a_{2j-1}$, Corollary \ref{coro:twice}
implies that, for all $1\leq i\leq n$, there is a a $\phi_M$-preserving isomorphism
$\psi_j:[0,b_j\setminus b_{j-1}]_{B(M)}\to[0,a_{2j}\setminus a_{2j-1}]_{B(M)}$.
Define $\psi(x)={\dot\bigvee}_{i=1}^n\psi_j(x\wedge(b_j\setminus b_{j-1}))$.

The proof that $\psi$ is a $\phi_M$-preserving isomorphism of Boolean algebras
is trivial and thus omitted.

\end{proof}
\begin{corollary}
\label{coro:firstiso}
Let $a,b\in B(M)$ be such that $\phi_M(a)=\phi_M(b)$. Then
there is a $\phi_M$-preserving 
isomorphism $\psi:[0,a]_{B(M)}\to [0,b]_{B(M)}$.
\end{corollary}
\begin{proof}
Use Lemma \ref{lemma:continuer} twice.
\end{proof}
\begin{lemma}
\label{lemma:basicauth}
Let $u,v\in B(M)$, $u\wedge v=0$, $\phi_M(u)=\phi_M(v)$.
Then there is a $\phi_M$-preserving automorphism $f$ of  $B(M)$
such that $f(u)=v$, $f(v)=u$ and for all
$x\leq (u\dvee v)^\complement$, $f(x)=x$.
\end{lemma}
\begin{proof}
By Corollary \ref{coro:firstiso}, there is an
isomorphism $\psi:[0,u]_{B(M)}\to [0,v]_{B(M)}$.
Let $f:B(M)\to B(M)$ be a mapping given by 
$$
f(x)=\psi^{-1}(x\wedge v)\dot\vee\psi(x\wedge u)
	\dot\vee(x\wedge(u\dot\vee v)^\complement).
$$

It is easy to check that, for all $x\in B(M)$, $f(f(x))=x$. Thus, $f$ is a bijection.
Moreover, we see that $f(0)=0$, $f(1)=1$ and, for
all $x,y\in B(M)$,
\begin{align*}
f(x\vee y)&=
   	\psi^{-1}\bigl((x\vee y)\wedge v\bigr)\dot\vee\psi\bigl((x\vee
	y)\wedge u\bigr)
	\dot\vee\bigl((x\vee y)\wedge(u\dot\vee v)^\complement\bigr)=\\
   	&=\psi^{-1}\bigl((x\wedge v)\vee (y\wedge
	v)\bigr)\dot\vee\psi\bigl((x\wedge u)\vee (y\wedge u)\bigr)
	\dot\vee\bigl(\bigl(x\wedge(u\dot\vee v)^\complement\bigr)\vee \\
	&\qquad \vee\bigl(y\wedge(u\dot\vee v)^\complement\bigr)\bigr)=\\
	&=\bigl(\psi^{-1}(x\wedge v)\dot\vee\psi(x\wedge u)
	\dot\vee(x\wedge(u\dot\vee v)^\complement)\bigr)\vee\\
	&\qquad\vee\bigl(\psi^{-1}(x\wedge v)\dot\vee\psi(x\wedge u)
	\dot\vee(x\wedge(u\dot\vee v)^\complement)\bigr)=\\
	&=f(x)\vee f(y)
\end{align*}
and
\begin{align*}
f(x^\complement)&=\psi^{-1}(x^\complement\wedge v)\dvee\psi(x^\complement\wedge u)
			\dvee\bigl(x^\complement\wedge(u \dvee v)^\complement\bigr)=\\
		&=\psi^{-1}\bigl(v\setminus(x\wedge v)\bigl)\dvee\psi\bigl(u\setminus(x\wedge u)\bigr)
			\dvee\bigl(x^\complement\wedge(u \dvee v)^\complement\bigr)=\\
		&=\bigl(u\setminus\psi^{-1}(x\wedge v)\bigr)\dvee\bigl(v\setminus\psi(x\wedge u)\bigr)
			\dvee\bigl(x^\complement\wedge(u \dvee v)^\complement\bigr)=\\
		&=\bigl(\psi^{-1}(x\wedge v)\dvee\psi(x\wedge u)
			\dvee\bigl(x\wedge(u \dvee v)\bigr)\bigr)^\complement
\end{align*}
The latter equality follows by elementary Boolean calculus.
Since $f$ preserves $0,1,\vee$ and $~^\complement$, it is a homomorphism of Boolean algebras.
\end{proof}
\begin{lemma}
\label{lemma:existsauth}
Let $u,v\in B(M)$, $\phi_M(u)=\phi_M(v)$.
Then there is a $\phi_M$-preserving automorphism $f$ of  $B(M)$
such that $f(u)=v$, $f(v)=u$ and for all
$x\leq (u\dvee v)^\complement$, $f(x)=x$.
\end{lemma}
\begin{proof}
Put $u_0=u\setminus u\wedge v$ and $v_0=v\setminus u\wedge v$.
Since 
$$
\phi_M(u_0)\oplus\phi_M(u\wedge v)=\phi_M(u)=\phi_M(v)=\phi_M(v_0)\oplus\phi_M(u\wedge v),
$$
$\phi_M(u_0)=\phi_M(v_0)$.
By Lemma \ref{lemma:basicauth}, there is
$f\in G(M)$ such that $f(u_0)=v_0$, $f(v_0)=u_0$ and for all
$x\in B$ such that $x\leq (u_0\dvee v_0)^\complement$ we have $f(x)=0$.
Since $u\wedge v\leq (u_0\dvee v_0)^\complement$, $f(u\wedge v)=u\wedge v$.
Therefore,
$$
f(u)=f\bigl(u_0\dvee(u\wedge v)\bigr)=f(u_0)\dvee(u\wedge v)=v_0\dvee(u\wedge v)=v
$$
and, similarly, $f(v)=u$.

Let $x\leq(u\vee v)^\complement$. Since
$x\leq(u_0\dvee v_0)^\complement$,
$f(x)=x$.
\end{proof}
\begin{corollary}
\label{coro:simGissimphi}
For all $u,v\in B(M)$, $u\sim_{G(M)} v$ iff
$\phi_M(u)=\phi_M(v)$.
\end{corollary}
\begin{proof}
One implication follows by the definition of $G(M)$, the other one
follows by Lemma \ref{lemma:existsauth}.
\end{proof}
\begin{corollary}
\label{coro:needed}
For all $u\in B(M)$, $u\sim_G \phi_M(u)$.
\end{corollary}
\begin{proof}
Put $v=\phi_M(u)$ in Corollary \ref{coro:simGissimphi}
\end{proof}

\begin{proof}[Proof of Theorem \ref{thm:MVtoMVP}]~

(MVP1):
Let $a,b\in B(M)$, $f\in G$ be such that $a\leq b$, $a\leq f(b)$.
Let $u=b\setminus(b\wedge f(b))$, $v=f(b)\setminus(b\wedge f(b))$.
We have
\begin{align*}
\phi_M(u)&=\phi_M\bigl(b\setminus(b\wedge f(b))\bigr)=\phi_M(b)\ominus\phi_M(b\wedge
		f(b)))=\\
         &=\phi_M(f(b))\ominus\phi_M(b\wedge f(b)))=
		\phi_M(f(b)\setminus(b\wedge f(b)))=\phi_M(v).
\end{align*}
By Lemma \ref{lemma:basicauth}, there is a $\phi_M$-preserving automorphism
$h$ of $B(M)$ with $h(u)=v$. Moreover, since $a\wedge u=a\wedge v=0$ and
$(b\wedge f(b))\wedge u=(b\wedge f(b))\wedge v=0$,
we have $h(a)=a$ and $h(b\wedge f(b))=b\wedge f(b)$.
This implies that
$$
h(b)=h((b\wedge f(b))\dot\vee u)=h((b\wedge f(b)))\dot\vee h(u)=
(b\wedge f(b))\dot\vee v=f(b).
$$
Thus, there is $h\in G$ such that 
$h(a)=a$ and $h(b)=f(b)$. By Lemma \ref{lemma:charmvp1}, this implies (MVP1).

(MVP2):
Let $a\wedge f(b)$ be an element of $L(a,b)$. By Corollary \ref{coro:needed},
there is $f_1\in G$ such that $f_1(a)=\phi_M(a)$.
Since $f_1$ is $\phi_M$-preserving, 
$\phi_M(f_1(a\wedge f(b)))=\phi_M(a\wedge f(b))$.
By Corollary \ref{coro:needed}, there is
$g\in G$ such that $g(f_1(a\wedge f(b)))=\phi_M(a\wedge f(b))$.
Since 
$$
f_1(a\wedge f(b))\leq f_1(a)=\phi_M(a)
$$
and
$$
g(f_1(a\wedge f(b)))=\phi_M(a\wedge f(b))\leq \phi_M(a),
$$
(MVP1) implies that there is $h\in G$ such that
$h(f_1(a\wedge f(b)))=\phi_M(a\wedge f(b))$ and $h(\phi_M(a))=\phi_M(a)$.

Put $y=a\wedge f_1^{-1}(h^{-1}(\phi_M(f(b))))$. We shall prove that
$y\geq a\wedge f(b)$ and that $y$ is a maximal element of $L(a,b)$.

Indeed, we have
$$
h(f_1(a))=h(\phi_M(a))=\phi_M(a),
$$
therefore
\begin{align*}
h(f_1(y))&=h\bigl(f_1\bigl(a\wedge f_1^{-1}(h^{-1}(\phi_M(f(b))))\bigr)\bigr)=\\
	&=h(f_1(a))\wedge h\bigl(f_1\bigl(f_1^{-1}(h^{-1}(\phi_M(f(b))))\bigr)\bigr)=\\
	&=\phi_M(a)\wedge \phi_M(f(b))=\phi_M(a)\wedge\phi_M(b)
\end{align*}
and
$$
h(f_1(a\wedge f(b)))=\phi_M(a\wedge(f(b)))\leq\phi_M(a)
	\wedge\phi_M(f(b))=h(f_1(y)).
$$
Since both $h$ and $f_1$ are automorphisms of $B(M)$, the latter inequality
clearly implies that $a\wedge f(b)\leq y$.
Moreover, since $h$ and $f_1$ are $\phi_M$-preserving and $\phi_M$ restricted
to $M$ is the identity mapping,
we obtain
$$
\phi_M(y)=\phi_M(h(f_1(y)))=\phi_M(\phi_M(a)\wedge\phi_M(b))
	=\phi_M(a)\wedge\phi_M(b).
$$

Let us prove that $y$ is maximal in $L(a,b)$.
Suppose that $z\in L(a,b)$, $z\geq y$.
Since $z=a\wedge f_2(b)$ for some $f_2\in G$, we see that
$$
\phi_M(z)=\phi_M(a\wedge f_2(b))\leq \phi_M(a)\wedge\phi_M(f_2(b))=\phi_M(y).
$$
This implies that $\phi_M(z)=\phi_M(y)$. As $\phi_M(z\setminus y)=
\phi_M(z)\ominus\phi_M(y)=0$ and $\phi_M$ is faithful,
$z\setminus y=0$ and hence $z=y$.

Let us prove that $\Alg(B(M),G(M))$ is isomorphic to $M$.
The isomorphism $\psi:\Alg(B(M),G(M))\to M$ is given by
$$
\psi([a]_{G(M)})=\phi_M(u).
$$
By Corollary \ref{coro:simGissimphi}, $\psi$ is well-defined and
injective. Since, for all $a\in M$,
$\psi([a]_{G(M)})=a$, $\psi$ is surjective. Obviously,
$\psi([1]_{G(M)})=1$. 
Let $[a]_{G(M)},[b]_{G(M)}\in\Alg(B(M),G(M))$ be such that
$[a]_{G(M)},[b]_{G(M)}$. We may always select the elements
$a,b\in B(M)$ so that $a\dvee b$ exists, that means, $a\wedge b=0$.
Since $\phi_M$ is a morphism of effect algebras, $\phi_M(a)\oplus\phi_M(b)$
exists in $M$ and we may compute
\begin{align*}
\psi([a]_{G(M)}\oplus[b]_{G(M)}&=\psi([a\dvee b]_{G(M)})=\phi_M(a\dvee b)=\\
	&=\phi_M(a)\oplus\phi_M(b)=\psi([a]_{G(M)})\oplus\psi([b]_{G(M)}),\\
\end{align*}
hence $\psi$ is a morphism of effect algebras. It remains to prove that $\psi$
is a full morphism. Suppose that $\psi([a]_{G(M)})\oplus\psi([b]_{G(M)})$
exists in $M$. Consider the elements $\phi_M(a)$ and
$\bigl(\phi_M(a)\oplus\phi_M(b)\bigr)\setminus\phi_M(a)$ of $B(M)$.  We see that
$$
\phi_M(a)\wedge\bigl(\bigl(\phi_M(a)\oplus\phi_M(b)\bigr)\setminus\phi_M(a)\bigr)=0,
$$
that means,
$\phi_M(a)\dvee\bigl((\phi_M(a)\oplus\phi_M(b)\bigr)\setminus\phi_M(a))$ exists in
$B(M)$. This implies that
$[\phi_M(a)]_{G(M)}\oplus[(\phi_M(a)\oplus\phi_M(b))\setminus\phi_M(a))]_{G(M)}$ exists in
$\Alg(B(M),G(M))$. Finally,
$$
\psi([\phi_M(a)]_{G(M)})=\phi_M(\phi_M(a))=\phi_M(a)=\psi([a]_{G(M)})
$$
and
\begin{multline*}
\psi([\bigl(\phi_M(a)\oplus\phi_M(b)\bigr)\setminus\phi_M(a)]_{G(M)})=
\phi_M\bigl(\bigl(\phi_M(a)\oplus\phi_M(b)\bigr)\setminus\phi_M(a)\bigr)=\\
=\phi_M\bigl(\phi_M(a)\oplus\phi_M(b)\bigr)\ominus\phi_M(\phi_M(a))
=\bigl(\phi_M(a)\oplus\phi_M(b)\bigr)\ominus\phi_M(a)=\\
=\phi_M(b)=\psi([b]_{G(M)}).
\end{multline*}
\end{proof}

\end{document}